\documentclass[reqno,11pt]{amsart}
\usepackage{amscd,amssymb}

\usepackage[all]{xy}

\CompileMatrices

\usepackage[dvipsnames]{xcolor}
\usepackage{graphics}
\usepackage{graphicx}

\theoremstyle{plain}
\newtheorem{Thm}[subsection]{Theorem}
\newtheorem{Cor}[subsection]{Corollary}
\newtheorem{Lem}[subsection]{Lemma}
\newtheorem{Prop}[subsection]{Proposition}
\newtheorem{Conj}[subsection]{Conjecture}

\theoremstyle{definition}
\newtheorem{Def}[subsection]{Definition}

\theoremstyle{remark}

\newtheorem{Rem}[subsection]{Remark}

\newtheorem{conj}{Conjecture}

\numberwithin{equation}{section}
\renewcommand{\rm}{\normalshape}

\newif\ifShowLabels
\ShowLabelstrue
\newdimen\theight
\def\TeXref#1{%
    \leavevmode\vadjust{\setbox0=\hbox{{\tt
        \quad\quad  {\small \rm #1}}}%
    \theight=\ht0
    \advance\theight by \lineskip
    \kern -\theight \vbox to
    \theight{\rightline{\rlap{\box0}}%
    \vss}%
    }}%

\ShowLabelsfalse

\renewcommand{\sec}[2]{\section{#2}\label{S:#1}%
    \ifShowLabels \TeXref{{S:#1}} \fi}
\newcommand{\ssec}[2]{\subsection{#2}\label{SS:#1}%
    \ifShowLabels \TeXref{{SS:#1}} \fi}

\newcommand{\refs}[1]{Section ~\ref{S:#1}}
\newcommand{\refss}[1]{Section ~\ref{SS:#1}}

\newcommand{\reft}[1]{Theorem ~\ref{T:#1}}
\newcommand{\refl}[1]{Lemma ~\ref{L:#1}}
\newcommand{\refp}[1]{Proposition ~\ref{P:#1}}

\newcommand{\refe}[1]{\eqref{E:#1}}

\newenvironment{thm}[1]%
    { \begin{Thm} \label{T:#1}  \ifShowLabels \TeXref{T:#1} \fi }%
    { \end{Thm} }

\renewcommand{\th}[1]{\begin{thm}{#1} \sl }
\renewcommand{\eth}{\end{thm} }

\newenvironment{lemma}[1]%
    { \begin{Lem} \label{L:#1}  \ifShowLabels \TeXref{L:#1} \fi }%
    { \end{Lem} }
\newcommand{\lem}[1]{\begin{lemma}{#1} \sl}
\newcommand{\elem}{\end{lemma}}

\newenvironment{propos}[1]%
    { \begin{Prop} \label{P:#1}  \ifShowLabels \TeXref{P:#1} \fi }%
    { \end{Prop} }
\newcommand{\prop}[1]{\begin{propos}{#1}\sl }
\newcommand{\eprop}{\end{propos}}

\newenvironment{corol}[1]%
    { \begin{Cor} \label{C:#1}  \ifShowLabels \TeXref{C:#1} \fi }%
    { \end{Cor} }
\newcommand{\cor}[1]{\begin{corol}{#1} \sl }
\newcommand{\ecor}{\end{corol}}

\newenvironment{defeni}[1]%
    { \begin{Def} \label{D:#1}  \ifShowLabels \TeXref{D:#1} \fi }%
    { \end{Def} }
\newcommand{\defe}[1]{\begin{defeni}{#1} \sl }
\newcommand{\edefe}{\end{defeni}}

\newenvironment{remark}[1]%
    { \begin{Rem} \label{R:#1}  \ifShowLabels \TeXref{R:#1} \fi }%
    { \end{Rem} }
\newcommand{\rem}[1]{\begin{remark}{#1}}
\newcommand{\erem}{\end{remark}}

\newenvironment{conjec}[1]%
    { \begin{Conj} \label{Co:#1}  \ifShowLabels \TeXref{Co:#1} \fi }%
    { \end{Conj} }
\renewcommand{\conj}[1]{\begin{conjec}{#1} \sl }
\newcommand{\econj}{\end{conjec}}

\newcommand{\eq}[1]%
    { \ifShowLabels \TeXref{E:#1} \fi
       \begin{equation} \label{E:#1} }
\newcommand{\eeq}{ \end{equation} }

\newcommand{\prf}{ \begin{proof} }
\newcommand{\epr}{ \end{proof} }

\newcommand\alp{\alpha}		
		
\newcommand\gam{\gamma}		
\newcommand\del{\delta}		\newcommand\Del{\Delta}

\newcommand\lam{\lambda}		\newcommand\Lam{\Lambda}

\newcommand\ome{\omega}		

\newcommand\calF{{\mathcal{F}}}

\newcommand\calH{{\mathcal{H}}}

\newcommand\calK{{\mathcal{K}}}
\newcommand\calL{{\mathcal{L}}}

\newcommand\calO{{\mathcal{O}}}

\newcommand\calS{{\mathcal{S}}}

		\newcommand\bfB{{\mathbf B}}
\newcommand\bfc{{\mathbf c}}

		\newcommand\bfG{{\mathbf G}}

		\newcommand\bfT{{\mathbf T}}
		\newcommand\bfU{{\mathbf U}}

		\newcommand\bfX{{\mathbf X}}

\newcommand\FF{\mathbb{F}}
\newcommand\GG{\mathbb{G}}

\newcommand\ZZ{\mathbb{Z}}

\newcommand\CC{\mathbb{C}}

	\newcommand\grb{{\mathfrak{b}}}

	\newcommand\grg{{\mathfrak{g}}}

	\newcommand\grn{{\mathfrak{n}}}

	\newcommand\grt{{\mathfrak{t}}}

\newcommand\sdp{\times \hskip -0.3em {\raise 0.3ex
\hbox{$\scriptscriptstyle |$}}} 

\newcommand\Aut{\operatorname{Aut}}

\newcommand\id{\operatorname{id}}

\newcommand\ord{\operatorname{ord}}

\newcommand\hatU{{\widehat{U}}}

\newcommand\x{\times}
\newcommand\ten{\otimes}

\newcommand{\ra}{\rangle}
\newcommand{\la}{\langle}

\newcommand\wt{\widetilde}

\newcommand{\kk}{\textsf{k}}

\renewcommand\kk{\Bbbk}
\newcommand\dep{\operatorname{dep}}
\newcommand{\st}{{\mathsf{t}}}
\renewcommand\O{\calO}
\newcommand\Heis{\operatorname{Heis}}
\newcommand{\iwk}{\mathsf{Iw}_K}
\newcommand{\hite}{\operatorname{ht}}
\newcommand\mc{\mathcal}
\newcommand{\md}{\operatorname{mod}}
\newcommand{\sph}{\operatorname{sph}}
\newcommand\mf{\mathfrak}
\renewcommand\ord{\operatorname{ord}}
\begin{document}

Dedicated to Igor Frenkel on the occasion of his 60th birthday

\bigskip

\title{An affine Gindikin-Karpelevich formula}
\author{A.~Braverman, H.~Garland, D.~Kazhdan and M.~Patnaik}
\maketitle

\begin{abstract}
In this paper we give an elementary proof of certain finiteness results about affine Kac-Moody groups over
a local non-archimedian field $\calK$. Our results imply those proven earlier in \cite{bk},\cite{bfk} and \cite{GaRo}
using either algebraic geometry or a Kac-Moody version of the Bruhat-Tits building.

The above finiteness results allow one to formulate an affine version of the {\em Gindikin-Karpelevich formula},
which coincides with the one discussed in \cite{bfk} in the case when $\calK$ has positive characteristic.
We deduce this formula from an affine version of the
 {\em Macdonald formula for the spherical function},
which will be proved in a subsequent publication.

\end{abstract}
\maketitle
\sec{int}{Introduction}
\ssec{}{Notations}
Let $\calK$ denote a local non-archimedian field with ring of integers $\calO$ and residue field $\kk$. Let
$\pi$ denote a generator of the maximal ideal of $\calO$.

Usually we shall denote algebraic varieties over $\calK$ (or a subring of $\calK$) by boldface letters $\bfX, \bfG$ etc.; their
sets of $\calK$-points will then be denoted $X,G$ etc.

Let $\bfG$ be a split, semi-simple, and simply connected algebraic group (defined over $\ZZ$) and let $\grg$ be its Lie
algebra.
\footnote{The reader should be warned from the very beginning that later we are going to change this notation; in particular,
later in the paper $\bfG$ will denote the corresponding affine Kac-Moody group}
As agreed above, we set $G=\bfG(\calK)$. In this Subsection we recall the usual Gindikin-Karpelevich formula for
$G$. Let $\bfT\subset \bfG$ be a maximal split torus; we denote its character lattice by $\Lam$ and its
cocharacter lattice by $\Lam^{\vee}$; note that since we have assumed that
$\bfG$ is simply connected, $\Lam^{\vee}$ is also the coroot lattice of $\bfG$. Given an element $x\in \calK^*$
we set $x^{\lam^{\vee}}=\lam^{\vee}(x)\in T$.

Let us choose a pair $\bfB,\bfB^-$ of opposite Borel subgroups such that
$\bfB\cap\bfB^-=\bfT$. We denote by $\bfU,\bfU^-$ their unipotent radicals. We denote by $R$ the set of roots of $\bfG$ and by
$R^{\vee}$ the set of coroots. Similarly $R_+$ (resp. $R^{\vee}_+$) will denote the set
of positive roots (resp. of positive coroots). We shall also denote by $2\rho$ (resp. by
$2\rho^{\vee}$) the sum of all positive roots (resp. of all positive coroots).
For any $\gam\in\Lam^{\vee}$ set $|\lam^{\vee}|=\la\lam^{\vee},\rho\ra$. Note that
if $\lam^{\vee}=\sum n_i\alp^{\vee}_i$ where $\alp^{\vee}_i$ are simple coroots, then
$|\lam^{\vee}|=\sum n_i$.
We let $W$ denote the Weyl group of $G$.

In addition we let $K=\bfG(\calO)$. This is a maximal compact subgroup of $G$.

For a finite set $X$, we shall denote by $|X|$ the number of elements in $X$.
\ssec{}{Gindikin-Karpelevich formula}
It is well-known that
$$
G=\bigsqcup\limits_{{\lam^{\vee}}\in\Lam^{\vee}_+} K\pi^{\lam^{\vee}} K.
$$

For any $\lam^{\vee}\in \Lam^{\vee}$ we set
\eq{c-function}
\bfc_{\grg} = \sum_{\mu^{\vee} \in \Lambda^{\vee}} | K \setminus K \pi^{\lambda^{\vee}} U^- \cap K \pi^{\lambda^{\vee} - \mu^{\vee}} U | e^{ \lambda^{\vee} - \mu^{\vee}} q ^{| \lambda^{\vee} - \mu^{\vee} | }.
\end{equation}
It is easy to see that the right hand side is independent of $\lam^{\vee}$; usually we shall take $\lam^{\vee}=0$.
\th{gk-finite}[Gindikin-Karpelevich formula]
\eq{gk-finite}
\bfc_{\grg}=\prod\limits_{\alp\in R_+}\frac{1-q^{-1}e^{-\alp^{\vee}}}{1-e^{-\alp^{\vee}}}.
\end{equation}
\eth

Let us note that classically, the Gindikin-Karpelevich formula computes the result of application of a certain intertwining operator
to the spherical vector in a principal series representation of $G$. However, it is trivial to see that the
classical Gindikin-Karpelevich formula is equivalent to the one presented above.
\ssec{}{Macdonald formula}
The proof of \refe{gk-finite} is actually very simple: it is easy to reduce it to the case
$\bfG=SL(2)$ where everything can be computed explicitly. However, the main purpose of this
paper is to present an analog of \refe{gk-finite} in the case when $\bfG$ is replaced
by an affine Kac-Moody group and in that case the above strategy does not work. So, let
us give a sketch of another proof of \refe{gk-finite}, based on the so called {\em Macdonald formula
for the spherical function}.

For any $\lam^{\vee}\in \Lam^{\vee}_+$ let us set
\eq{formal:sph}
\calS(\lambda^{\vee}):= \sum_{\mu^\vee \in \Lambda^{\vee} } | K \setminus K \pi^{\lambda^\vee} K \cap K \pi^{ \mu^\vee} U | q^{ \langle \rho, \mu^{\vee} \rangle} e^{\mu^{\vee}}.
\end{equation}
This is an element of the group algebra $\CC[\Lam^{\vee}]$. Actually, general facts about the Satake isomorphism (see \refss{satake})
imply that $\calS(\lam^{\vee})\in\CC[\Lam^{\vee}]^W$.

\th{macdonald-finite}[Macdonald]
Let us set
$$
\Delta = \prod_{\alpha\in R_+} \frac{1 - q^{-1} e^{ - \alpha^\vee} } {1 - e^{ - \alpha^\vee}}.
$$
Then for any $\lambda^\vee\in \Lam^{\vee}_+$ we have
\eq{sph:mac}
\calS(\lambda^{\vee}) =
\frac{q^{ \langle \rho, \lambda^{\vee} \rangle}}{W_{\lam^{\vee}}(q^{-1}) }\sum_{w \in W} w (\Delta) e^{w \lambda^{\vee}}
\end{equation}
Here $W_{\lam^{\vee}}$ is the stabilizer of $\lam^{\vee}$ in $W$ and
$$
W_{\lam^{\vee}}(q^{-1})=\sum \limits_{w\in W_{\lam^{\vee}}} q^{-\ell(w)}.
$$
\eth

It is easy to see that Macdonald formula implies \refe{gk-finite}. The reason is that (as is well-known)
for any fixed $\mu^\vee$ we have
\eq{bg}
K \pi^{\lambda^\vee} K \cap K \pi^{\lambda^\vee - \mu^{\vee} } U = K \pi^{\lambda^{\vee}} U^- \cap K \pi^{\lambda^{\vee} - \mu^\vee} U \end{equation}
for any  sufficiently dominant $\lambda^{\vee}$.
Using \refe{bg} to take the limit of \refe{formal:sph} as $\lambda^{\vee}$ is made more and more dominant (which we write as $\lambda^{\vee} \to + \infty$) one gets \refe{gk-finite}.
\ssec{satake}{Interpretation via Satake isomorphism}
Let us explain the meaning of \refe{formal:sph} in terms of the {\em Satake isomorphism}.
Let $\calH_{\sph}$ denote the spherical Hecke algebra of $G$. By definition this is the convolution algebra
of $K$-bi-invariant distributions with compact support on $G$; by choosing a Haar measure on $G$ for which
the volume is $K$ is equal to 1, we may identify it with the space of $K$-bi-invariant functions with
compact support. Thus the algebra $\calH_{\sph}$ has a basis $h_{\lam^{\vee}}$ where
$h_{\lam^{\vee}}$ is the characteristic function of the corresponding double coset. The Satake isomorphism is
a canonical isomorphism  between $\calH$ and $\CC[\Lam^{\vee}]^W$. By its very definition
this isomorphism sends $h_{\lam^{\vee}}$ to $\calS(\lam^{\vee})$. Thus \refe{sph:mac} can be thought of
as an explicit computation of the Satake isomorphism in terms of the above basis. In particular, since
$h_0=1$, the right hand side of \refe{sph:mac} for $\lam^{\vee}=0$ must be equal to 1; this is not
completely obvious from the definition.

\ssec{}{Affine Kac-Moody groups}
Our aim in this paper is to discuss a generalization of the Gindkin-Karpelevich formula to the
case when $\bfG$ is an untwisted affine Kac-Moody group. When $\calK$ has positive characteristic
this was already done in \cite{bfk} (modulo the assumption that the results of \cite{bfg} hold in positive characteristic).

Let $\bfG_o$ be a simple, simply connected group. Then one can form the polynomial loop group $\bfG_o[t, t^{-1}]$ which admits a central extension by $\mathbb{G}_m$ which we denote by $\wt{\bfG}.$ The full affine Kac-Moody group is then $\bfG =\wt{\bfG} \rtimes \mathbb{G}_m$ where $\mathbb{G}_m$ acts by rescaling the loop parameter $t.$ We shall be interested in the group $G=\bfG(\calK)$. In this setting, we may define analogues of the groups $K, U, U^-$.

\th{finiteness}
With the above notation we have:
\begin{enumerate}
  \item
  For $\lambda^{\vee}, \mu^\vee \in \Lambda^{\vee}$ the set
  $K \setminus K \pi^{\lambda^\vee} U^-\cap K \pi^{\mu^\vee} U$ is finite.
  Moreover, it is empty unless $\lam^{\vee}\geq \mu^{\vee}$.
  \item
  For any $\lam^{\vee}\in \Lam^{\vee}_+$, $\mu^{\vee}\Lam^{\vee}$ as above the set
$$
K\backslash K \pi^{\lambda^\vee} K \cap K \pi^{ \mu^{\vee} } U
$$
  is finite. Moreover, it is empty unless $\lam^{\vee}\geq \mu^{\vee}$.

\item The set $K \pi^{\lambda^\vee} K \cap K \pi^{ \mu^{\vee} } U^-$ is empty unless
$\lam^{\vee}\geq \mu^{\vee}$.
\item
  For $\lambda^{\vee}$ sufficiently dominant and fixed $\mu^{\vee}$ we have
  \eq{bg-aff}
  K \pi^{\lambda^\vee} K \cap K \pi^{\lambda^\vee - \mu^{\vee} } U = K \pi^{\lambda^{\vee}} U^- \cap K \pi^{\lambda^{\vee} - \mu^\vee} U
  \end{equation}
\end{enumerate}
\eth
The third statement of \reft{finiteness} is easy; it is also easy to see that (1) and (3) imply (2).
Statement (4) is also not difficult, once we know (1) (details
can be found in \refs{final}). Thus essentially it is enough to prove (1).
This is the main result of this paper.
To prove it, we develop some techniques for analyzing the infinite dimensional group $U^-: $ the main idea is to break up the group $U^-$ into pieces $U^-_w$ indexed by the Weyl group of $\bfG$. These pieces are not finite-dimensional, but if we replace $U^-$ in the above intersection by $U^-_w$ the resulting set is finite. We can then also show, using a certain result from representation theory, that only finitely many $w$ contribute to the above intersection.

We should note right away that \reft{finiteness}(2) appears in \cite{bk}. In addition, while this paper was in preparation, the paper \cite{GaRo} has appeared, where an analog of \reft{finiteness}(2)
is proved in a greater generality (for $\bfG$ being any almost split symmetrizable Kac-Moody group).
However, the proofs in both \cite{bk} and \cite{GaRo} are quite cumbersome and use some rather complicated machinery.
From the results of this paper we get a new elementary (though not very easy) proof of \reft{finiteness}(2).
Let us also add that if one assumes \reft{finiteness}(4) that clearly \reft{finiteness}(1) follows
from (the earlier known) \reft{finiteness}(2). However, we don't know how to prove \reft{finiteness}(4) without
using \reft{finiteness}(1).

Using the above theorem, we can then make sense in our infinite dimensional setting of the formal Gindikin-Karpelevich sum $\mathbf{c}_{\grg}.$ Let us now formulate the statement of the affine Gindikin-Karpelevich formula. For this we first want
to discuss the affine Macdonald formula.
\ssec{}{Affine Macdonald formula}
Let us define $\calS(\lam^{\vee})$ as in \refe{formal:sph}. Note that it makes sense due to \reft{finiteness}(2).

Let us also set
\eq{hlam}
H_{\lam^{\vee}}=\frac{\st^{ \langle \rho, \lambda^{\vee} \rangle}}{W_{\lam^{\vee}}(\st^{-1}) }\sum_{w \in W} w (\Delta) e^{w \lambda^{\vee}}.
\end{equation}
Here $\st$ is a formal variable (which should not be confused with $t$!),  $\Del$ is defined as in the finite-dimensional case (with $q$ replaced by $\st$) and $\rho$ is the element of $\Lam$ such that $\la \rho,\alp_i^{\vee}\ra =1$
for any simple coroot $\alp_i^{\vee}$.
\th{aff-mac}
$\calS(\lam^{\vee})$  is equal to the specialization of $\frac{H_{\lam}}{H_0}$ at ${\st=q}$.
\eth
This theorem will be proved in \cite{bkp}.

Note that contrary to the finite dimensional case, $H_0\neq 1$. In fact, the function $H_0$ was studied by Macdonald in
\cite{ma:formal} using the works of Cherednik. Macdonald has shown that $H_0$ (which a priori is defined as an infinite sum)
has an infinite product decomposition. For example, when $\grg_o$ is simply laced, Macdonald's formula reads as follows:
\eq{mac2003}
H_0=\prod\limits_{i=1}^{\ell}\prod\limits_{j=1}^{\infty} \frac{1-\st^{-m_i}e^{-j\del}}{1-\st^{-m_i-1}e^{-j\del}}.
\end{equation}
Here $\del$ is the minimal positive imaginary coroot of $\grg$ and $m_1,\cdots,m_{\ell}$ are the exponents of $\grg_o$.

A similar product decomposition for $H_0$ exists for any $\grg_o$ (cf. \cite{ma:formal}).
\ssec{}{Affine Gindikin-Karplevich formula}
We are now ready to formulate the affine version of the Gindikin-Karpelevich formula.
\th{gk-aff}
For any $\grg$ as above, we have
\eq{gk-aff}
\bfc_{\grg}=\frac{1}{H_0}\prod\limits_{\alp\in R_+}\Bigg(\frac{1-q^{-1}e^{-\alp^{\vee}}}{1-e^{-\alp^{\vee}}}\Bigg) ^{m_{\alp}}.
\end{equation}
Here the product is taken over all positive roots (or, equivalently, coroots) of $\grg$ and $m_{\alp}$ denotes the multiplicity
of the coroot $\alp^{\vee}$.
\eth

Let us stress that \reft{gk-aff} depends on \reft{aff-mac} which will be proved in another publication.
Also, when $\calK=\FF_q((t))$, \reft{gk-aff} was proved in \cite{bfk} modulo an assumption that certain geometric
results from \cite{bfg} hold in positive characteristic (\cite{bfg} deals only with characteristic 0).

We conclude with a final remark regarding the constant term conjecture of Macdonald (Cherednik’s Theorem),
which is essentially the statement that $H_0^{-1}$
is the imaginary root part (or “constant term” in the terminology of \cite{mac:mad}, \cite{mac:cox}
of $\Del^{-1}$: using the techniques of \cite{bkp} and this paper, we can show (without resorting to the work of Cherednik)
that the constant term of
$\Del^{-1}$ can be related to the Gindikin-Karpelevic sum $\bfc_{\grg}$. On the other hand, in the case
when the local field $\calK$ has positive characteristic, the Gindikin-Karpelevic sum $\bfc_{\grg}$ has been computed in \cite{bfk}
using a geometric argument which is independent of Cherednik’s work. Thus the combination of the works of
\cite{bkp,bfk} and this paper also yield an independent proof of the constant term conjecture of Macdonald.

\ssec{}{Organization of the paper}In \refs{aff-KM} we recall some basic facts about affine Kac-Moody group schemes.
\refs{part1}, \refs{completion} and \refs{part2} are devoted to the proof of the first assertion of \reft{finiteness}.
In \refs{final} we prove the other assertions of \reft{finiteness} and explain how to deduce \reft{gk-aff} from \reft{finiteness}
and \reft{aff-mac}.
\ssec{}{Acknowledgements}A.~B. was partially supported by the NSF grant DMS-0901274.
A.B. and D.K. were partially supported by a BSF grant 2008386.
M.P. was supported by an NSF Postdoctoral Fellowship, DMS-0802940 while this work was
being completed.
\sec{aff-KM}{Affine Kac-Moody groups}

\ssec{}{Affine Kac-Moody group functor}
Let $\bfG_o$ be a split simple simply connected group over $\ZZ$; let $\grg_o$ denote its Lie algebra. To this group one can attach a group ind scheme
$\bfG$ (called the affinization of $\bfG_o$) in the following way.

First, one considers the formal loop group functor $\bfG_o((t))$ (by definition it sends a commutative
ring $R$ to $\bfG_o(R((t)))$). Clearly, the multiplicative group $\GG_m$ acts on $\bfG_o((t))$ and we denote
by $\bfG'$ the corresponding semi-direct product. It is well-known (cf. e.g. \cite{fal}) that any invariant
integral bilinear form on $\grg_o$ gives rise to a central extension of $\bfG'$ by means of $\GG_m$.
In the case, when the bilinear form in question is the minimal one (equal to the Killing form
divided by $2h^{\vee}$ where $h^{\vee}$ is the dual Coxeter number of $\bfG_o$) we denote the corresponding
central extension by $\bfG$. Thus, by definition we have a short exact sequence of group ind-schemes
$$
1\to \GG_m\to \bfG\to \bfG'\to 1.
$$

The above central extension is known to be split over $\bfG_o[[t]]$ and in what follows we choose such a splitting.

We choose a pair of opposite Borel subgroups $\bfB_o,\bfB_o^-$ with unipotent radicals $\bfU_o,\bfU_o^-$.
The intersection $\bfT_o=\bfB_o\cap \bfB_o^-$ is a maximal torus of $\bfG_o$.

Let $\bfT=\GG_m\x \bfT_o\x\GG_m$.
Let $\bfG_o[[t]]_\bfB$ denote the preimage of $\bfB_o$ under the natural map $\bfG_o[[t]]\to\bfG_o$. We let
$\bfB$ to be the preimage in $\bfG$ of $\bfG_o[[t]]_{\bfB}\rtimes \GG_m\subset \bfG'$.
This is a group-scheme, which is endowed with a natural map to $\bfT$. We denote by $\bfU$ the kernel of this map.
This is the pro-unipotent radical of $\bfB$.

Similarly, let $\bfG_o[t^{-1}]_{\bfB^-}$ be the preimage of $\bfB_o^-$ under the map
$\bfG_o[t^{-1}]\to\bfG_o$ coming from evaluating $t$ to $\infty$. We let
$\bfB^-\subset \bfG$ to be the preimage of $\bfG_o[t^{-1}]_{\bfB^-}\rtimes \GG_m\subset \bfG'$.
This is a group ind-scheme, which (similarly to $\bfB$) is endowed with a natural map to $\bfT$ and
we denote its kernel by $\bfU^-$.

In addition, the intersection $\bfB\cap \bfB^-$ is naturally isomorphic to $\bfT$.

\ssec{}{The affine root system}
We denote the Lie algebras of $\bfG,\bfT,\bfB,\bfB^-,\bfU,\bfU^-$ respectively by
$\grg,\grt,\grb,\grb^-,\grn,\grn^-$.
We shall denote by $R_o$ the set of roots of $\grg_o$ and by $R$ the set of roots of $\grg$.
Although there is a natural embedding $R_o\subset R$, in the future it will be convenient to
use different notations for elements of $R_o$ and $R$: we shall typically denote elements of
$R_o$ $\alp,\beta,...$ and elements of $R$ by $a,b,...$.

Let $\Pi= \{ a_1, \ldots, a_{l+1} \}$ be the set of simple roots such that for $i\leq l$ the element $a_i$ is a simple root in $R_o$
(which we shall also denote by $\alp_i$).   We define the \emph{height} of a root $\alpha \in R_o$  as follows: if $\alpha = \sum_{i=1}^l p_i \alpha_i,$ then
$$
\hite(\alpha) = \sum_{i=1}^l p_i.
$$
Recall that $a_{l+1}= - \theta + \delta$ where $\theta$ is the classical root of maximal height and $\delta$ is the minimal positive imaginary root.   Let $Q$ be the (affine) root lattice. Any $\gam\in Q$ we may write it uniquely as
$$
\gam = \sum_{i=1}^{l+1} p_i a_i.
$$

We define the \emph{depth of $\gam$}  to be
$$
\dep(\gam) = -p_{l+1}.
$$
 If we have $\dep(\gam)=0$ (i.e. $p_{l+1}=0$), then we define the \emph{order of $\gam$} to be
$$
\ord(\gam) = -\hite(\gam).
$$
Note that if we write $a \in R$ as $a = \alpha + n \delta$ for $\alpha \in R_o,$ then $\dep(a)=-n.$

We denote by $R_{o,+}$ (resp. $R_{o,-}$ the set of positive (resp. negative) roots of $\grg_o$.
For each $\alp\in R_{o,+}$ (resp. $\alp\in R_{o,-}$ we denote by $\chi_{\alp}:\GG_a\to \bfU_o$ (resp.
$\chi_{\alp}:\GG_a\to \bfU_o^-$) the corresponding homomorphism. We shall use similar notation for
real roots of $\grg$.
\ssec{}{The Chevalley basis and the integral form}
Let us pick a Chevalley basis for $\grg$ which we write as follows: if $a= \alpha - n \delta \in R_{re, -}$ then the basis element in $\grg^{a}$ is denoted by $\xi_a= \xi_{-n, \alpha}.$ Let $\{ \xi_i \}_{i=1}^l$ denote a Chevalley basis for $\grt_o$ and
if $a= - n \delta$ then we denote by $\{ \xi_{n, i} \}_{i=1}^l$ the Chevalley basis for $\grg^{-n \delta}=\grt_o \otimes t^{-n}.$

This basis also defines an integral form $U_{\ZZ}(\grg)$ of the universal enveloping algebra of $\grg$ as in \cite{Kos}.
It has the triangular decomposition
$$
U_{\ZZ}(\grg)=U_{\ZZ}(\grn^-)\ten U_{\ZZ}(\grt)\ten U_{\ZZ}(\grn).
$$

\ssec{}{Integrable modules}
We let $\Lam$ denote the lattice of characters of $\bfT$; let also $\Lam^{\vee}$ denote the dual lattice.
Let $\Lam^+$ denote the set of dominant weights (and $(\Lam^{\vee})^+$ the set of dominant coweights).
For any $\lam\in\Lam^+$ we let $V^\lam_{\ZZ}$ denote the corresponding integrable highest weight $\bfG(\ZZ)$-module (a.k.a. Weyl module).
It is defined over $\ZZ$ and we shall write $V^{\lam}_R$ for $V_{\ZZ}^{\lam}\ten R$ (this is a $\bfG(R)$-module).
We shall usually set $V^{\lam}=V^{\lam}_{\calK}$.

The module $V^{\lam}_{\ZZ}$ has a weight decomposition
$$
V^\lam_{\ZZ}=\oplus V^{\lam}_{\ZZ}(\mu).
$$
For every $w\in W$ the weight module $V^{\lam}_{\ZZ}(w\lam)$ is of rank 1. For every such $w$ we set
$V_{w,\ZZ}^{\lam}=U_{\ZZ}(\grb)\cdot V_{\ZZ}^{\lam}(w\lam)$. This is a free $\ZZ$-module of finite rank inside
$V^{\lam}_{\ZZ}$ (usually called "the Demazure module"). Hence $V^{\lam}_w$ is a finite-dimensional $
\calK$-vector space
(endowed with an action of $B$).

In addition, let us set $V^{\lam}(m)$ to be the direct sum of all the $V^{\lam}(\mu)$ for which $\dep(\mu)>m$.
Also we set $V^{\lam}[m]=V^{\lam}/V^{\lam}(m)$. This is a finite-dimensional $U^-$-submodule.

\ssec{}{The Heisenberg subalgebra}
The algebra $\grg$ also contains a rank $l$ infinite-dimensional \emph{Heisenberg subalgebra} $\Heis \subset \grg,$
$$
\Heis:= \mathcal{K} \textsf{c} \oplus \bigoplus_{n \in \ZZ} \grt_o \otimes t^n
$$
where the element $\textsf{c}$ is the central element in $\grt$.
We shall say that an element in $\grt_o \otimes t^n$ is of \emph{degree} $n$.
\ssec{}{Some decompositions}
Let now $\calK$ be a local non-archimedian field with ring of integers $\calO$ and residue field $\kk$.
We let $q$ denote the cardinality of $\kk$; we also choose a generator of the maximal ideal of $\calO$ which we
shall denote by $\pi$.

For a scheme $\bfX$ over $\calO$ we set $X=\bfX(\calK), X_{\calO}=\bfX(\calO), X_{\kk}=\bfX(\kk)$.

Let us set $I$ (resp. $I^-$) to be the preimage of $B_{\kk}$ (resp. $B^-_{\kk}$) under the natural map $G(\calO)\to G_{\kk}$.
We shall choose a lift of all the elements of $W$ to $K$ and we are going to call this set of representatives
$\dot{W}$. We shall sometimes identify $\dot{W}$ with $W$ (when it doesn't lead to a confusion).

We have the Bruhat decompositions
$$
\bfG_{\kk} = B_{\kk} \dot{W} B_{\kk}; \quad G_{\kk} = B_{\kk} \dot{W} B^-_{\kk}
$$
which lift to give decompositions of the form
$$
K = I \dot{W} I  \text{ and }  K = I^- \dot{W} I
$$
In addition we have the Iwasawa decomposition
$$
G=K\cdot B
$$
(cf. e.g. Section 2 of \cite{bfk} for the proof).


\sec{part1}{Finiteness of the Gindikin-Karpelevich Sum, Part I}

The main  goal of this paper will be to prove the following result which shows that the sets appearing in the expression \refe{c-function}
are finite.

\th{gk:fin} Let $\mu^{\vee} \in Q_+.$ Then we have,
\label{gk:fin}
\eq{gk:fin:1} | K \backslash K U^- \cap K \pi^{-\mu^{\vee}} U | <  \infty.
\end{equation}
\eth

\noindent The proof will be carried out in three steps:
\begin{enumerate}
  \item
  We shall decompose the set $U^-$ into disjoint subsets $U^-= \bigcup_{w \in W} U_w^-;$
  \item
  We shall show that there are only finitely many $w \in W$ such that $U^-_w$ can contribute to \refe{gk:fin:1}
  \item
  We shall show that the sets \refe{gk:fin:1} with $U^-$ replaced by $U^-_w$ are finite;
\end{enumerate}

Steps 1 and 2 will be performed in this Section. Step 3 will be performed in \refs{part2}.

\noindent
{\em Step 1:} We have an injection $U^{-} \hookrightarrow G/B.$ By the Iwasawa decomposition, we have that
$G/B \cong K / K \cap B$ and so we have an embedding,
$$
\iwk: U^- \hookrightarrow K / K \cap B.
$$
More explicitly, if we write $u^- = k h u$ in terms of its Iwasawa coordinates, then the above map sends $u^- \mapsto k\; (\text{mod}\; K \cap B).$

Since $K \setminus K U^- = U^-_{\calO} \setminus U^-,$ the map $\iwk$ descends to a map
$$
 \iwk: K \setminus K U^- \hookrightarrow U^-(\O) \setminus K / K \cap B.
$$
Let us denote by $\varpi$ the reduction modulo $\pi$ map $\varpi: K \to G_{\kk}$.

By composing the map
$\iwk:K \setminus K U^- = U^-_{\O} \setminus U^- \to U^-_{\O} \setminus K / K \cap B$ with the the natural map
$ U^-_{\O} \setminus K / K \cap B\to U^-_{\kk} \setminus G_{\kk} / B_{\kk} = W$ we get a map
$\varphi: K \setminus K U^- \to W$.
Letting $\psi: U^- \to U^-_{\O} \setminus U^-$ denote be the projection map, we define
$$
U^-_w =  \psi^{-1} \phi^{-1}(w) \subset U^-.
$$

We can give a more explicit description of the elements in $U^-_w$ as follows.  First set
$$
G_{\pi} = \{ k \in K | \varpi(k) = 1. \}
$$
We can define $U_{\pi}, U^{-}_{\pi}$ and $T_{\pi}$ in the same way.

Then $G_{\pi} \subset I$ and in fact we have a direct product decomposition
$G_{\pi} = U_{\pi} U^-_{\pi} T_{\pi}$.
We have
$$
\varpi^{-1}(w) = U^{-}_{\O} G_{\pi} w B_{\O} = U^-_{\O} U_{w, \pi} w B_{\O}
$$
where
$$
U_{w, \pi} = \prod_{a > 0, w^{-1} a < 0 } U_{a, \pi}
$$
where
$U_{a, \pi}$ is the first congruence subgroup in $U_{a,\calO}$.

So every $u^- \in U^-_w$ has an expression of the form
\eq{uw}
 u^- = u^-_{\O} u_{w, \pi} w b_{\O} \pi^{\xi} u
\end{equation}
where lowercase elements are in the corresponding uppercase subgroups and $\xi \in \Lambda^{\vee}$.
\lem{step1}
For $u_{w, \pi} \in U_{w, \pi},$ there exist $u^- \in U^-, u \in U, h \in T$ such that
$$
u_{w, \pi} w = u^- h u.
$$
Moreover, $u^-, u, h$ are uniquely defined.
\elem
\prf
If
\eq{uw:mu:2}
u_{w, \pi} = u_1^- h_1 u_1 = u_2^- h_2 u_2 \; \text{ where } u_1^{\pm}, u_2^{\pm} \in U^{\pm}, \; h_1, h_2 \in T
\end{equation}
then it follows that $(u_2^-)^{-1} u_1^- \in B,$ which is a contradiction unless $u_1^- =u_2^-.$
Similarly, we conclude that $u_1 = u_2$ and $h_1=h_2.$
\epr

\noindent
{\em Step 2:} We now prove the following inequality which is the main result of this Section.
\lem{drin}
Fix $\mu^{\vee} \in Q^{\vee}_+$ Then if
$$
K U^-_{w} \cap K \pi^{-\mu^{\vee}} U \neq \emptyset,
$$
we must have
$$
\frac{l(w)}{2} \leq |\mu^{\vee}|:= \la \rho, \mu^{\vee} \ra.
$$
\elem
\prf[Proof of Lemma]   Suppose $u^-_w \in U^-_w$ is such that $Ku^-_w \in KU^-_w \cap K \pi^{-\mu^{\vee}} U.$ Then we can write
$$
u^-_{w} = u^{-}_{\O} u_{w, \pi} w b_{\O}  \pi^{-\mu^{\vee}} u,
$$
which produces a relation of the form
\eq{step3:1} u^- = u_{w, \pi} w b_{\O}  \pi^{-\mu^{\vee}} u.
\end{equation}

Applying both sides of \refe{step3:1} to $v_{\rho},$ a primitive highest weight vector in  $V^{\rho}:$
\eq{step3:2}
u^- v_{\rho} = u^{+}_{w, \pi} w b_{\O} \pi^{-\mu^{\vee}} u v_{\rho}.
\end{equation}
The left hand side of \refe{step3:2} is of the form
\eq{step3:2.1}
v_{\rho} + \text{lower terms},
\end{equation}
whereas the right hand side is of the form
\eq{step3:2.2}
\delta \pi^{- \langle \rho, \mu^{\vee} \rangle} u_{w, \pi} v_{w \rho},
\end{equation}
where $\delta \in \O^*$ and $v_{w \rho}$ is again a primitive vector in $V_{\O}^\rho.$

Consider an element $u_{w, \pi} \in U_{w, \pi}$ and keep the notation of Step 2. So the element $u_{w, \pi}$ acts acts via a sum of the form,
\eq{step3:3}
\sum_{n_1, \ldots, n_r} \sigma_1^{n_1} \cdots \sigma_r^{n_r} \xi_{\beta_1}^{(n_1)} \cdots \xi_{\beta_r}^{(n_r)},
\end{equation}
where the $\beta_i,$ for $i=1, \ldots, r$ are the positive real roots $a \in R_{re, +}$ such that $w^{-1} a <0,$ and $\xi_{\beta_i}^{(n_i)}$ are the divided powers of our fixed Chevalley basis elements.
Let us now consider the element $v_{\rho}$ in the highest weight space of $V^{\rho}$
 in the expression \refe{step3:2.2}.

Let us now use the following result due to A.~Joseph (cf. \cite{bfg}, Lemma 18.2):
\prop{jos}
Suppose that $v_{\rho} \in F^{m}(U (\mathfrak{n}^+))v_{w \rho}.$ Then,
$$
m \geq l(w)/2,
$$
 where $l(w)$ is the length of $w.$
\eprop

By \refp{jos} we must apply at least $l(w)/2$ operators from $\mathfrak{n}$ to $v_{w \rho}$ in order to obtain an element in $V_{\rho}.$ This corresponds, in an expression of the form \refe{step3:3}, to terms in which $n_1 + \cdots + n_r= m \geq l(w)/2.$  Since the $\xi_{\beta_i}^{(n)}$ map $V^{\rho}(\calO)$ into itself, such an expression will introduce a zero of order at least $n_1 + \cdots + n_r$ into the resulting element in $V^\rho.$ However, since from \refe{step3:2.1} the element produced in $V^\rho$ as a result of applying $u_{w, \pi}$ must be primitive, we obtain the desired equality:
 $$
 |\mu^{\vee}| \geq l(w)/2
 $$
 in light of \refe{step3:2.2}.

\epr
Thus to prove \reft{gk-aff}(1) we need to show that the quotient $K\backslash KU_w^-\cap K\pi^{-\mu^{\vee}}U$ is finite
for every $w$ and $\mu^{\vee}$. In other words, we need to prove the following
\th{main-fin}
The set $K\backslash KU^-(w,\mu^{\vee})$ is finite, where
\eq{uwm}
U^-(w,\mu^{\vee})=U^-\cap U_{w,\pi}T_{\calO}w\pi^{-\mu^{\vee}}U.
\end{equation}
\eth
This will be done in \refs{part2}.
 But first we need to introduce some notation
related to the group $U^-$ and its completion.

\sec{completion}{Completion of $U^-$}

The purpose of this section is to discuss some sort of coordinates on a formal completion of $U^-$ which will be used later.
\ssec{}{The completion}
Recall that the group ind-scheme $\bfU^-$ is the the preimage of $\bfU^-_o$ under the natural (evaluation at $\infty$)
map $\bfG_o[t^{-1}]\to \bfG_o$. Let ${\mathbf \hatU}^-$ denote the preimage of $\bfU^-_o$ in $\bfG[[t^{-1}]]$. This is
a group-scheme; we have a natural map $\bfU^-\to{\mathbf \hatU}^-$ which
induces an injection
$$
i: U^- \hookrightarrow \hatU^-.
$$
We shall often identify $U^-$ with its image in $\hatU^-.$

For every $m\geq 0$ we shall set $U^-(m)$ to be the subgroup of $U^-$ consisting of elements which are equal to $1$
modulo $t^{-m}$. We set $U^-[m]=U^-/U^-(m)$ and we shall denote by $\ome_m$ the natural projection from
$U^-$ to $U^-[m]$.
\ssec{}{Some infinite products}
For $\alpha \in R_{o, +}$ and $\beta \in R_{o, -}$ let us set
$$
\sigma_{\alpha}= \sum_{i=1}^{\infty} c_i t^{-i} \in t^{-1}\mathcal{K}[[t^{-1}]]  \text{ and }  \sigma_{\beta}= \sum_{i=0}^{\infty} c_i t^{-i} \in \mathcal{K}[[t^{-1}]].
$$
Then we may consider the following  products as elements of $\hatU^-$
$$
\chi_{\alpha}(\sigma_{\alpha}) = \prod_{i=1}^{\infty} \chi_{\alpha - i \delta} (c_i)  \text{ and }  \chi_{\alpha}(\sigma_{\alpha}) = \prod_{i=0}^{\infty} \chi_{\beta - i \delta} (c_i).
$$
Suppose we are given a unit $\sigma \in \mathcal{K}[[t^{-1}]]^*.$ If $\sigma \equiv 1 (\md t^{-1}),$ we have a factorization,
\eq{fact} \sigma = \prod_{j \geq 1} (1+ c_jt^{-j}),
\end{equation}
where the $c_j$ are uniquely determined.  For $i=1, \ldots, l$ we form the expressions,
$$
h_i(\sigma) := \chi_{\alpha_i}(\sigma)\chi_{- \alpha_i} (-\sigma^{-1}) \chi_{\alpha_i} (\sigma) \chi_{\alpha_i}(1) \chi_{-\alpha_i}(-1) \chi_{\alpha_i}(1)
$$
which again define elements of $\hatU^-.$ With respect to the factorization \refe{fact} we then have
$$
h_i(\sigma) = \prod_{j \geq 1} h_i(1+ c_jt^{-j}).
$$

Fix a positive integer $m \geq 1 $ and consider now an element of $\hatU^-$ of the form
\eq{im} u^-[m]:= \prod_{\alpha \in R_{o, +} } \chi_{\alpha}(t^{-m} s_{m, \alpha} ) \; \prod_{i=1}^l h_i(1+c_{i, m}t^{-m}) \; \prod_{\alpha \in R_{o, +}} \chi_{- \alpha} (t^{-m} \tilde{s}_{m, \alpha} ),
\end{equation}
where the products are with respect to a fixed ordering
on $R_{o,+}$ and $s_{m, \alpha},\; \tilde{s}_{m, \alpha},\; c_{i, m} \in \mathcal{K}.$ If $m=0,$ consider elements of the form
\eq{im:0}
u^-[0] := \prod_{\alpha \in R_{o,+}} \chi_{- \alpha} (\tilde{s}_{0, \alpha}) \; \text{ with  } \tilde{s}_{0, \alpha} \in \mathcal{K}.
\end{equation}

For future use, we shall need to understand the action of the various components of $u^-[m]$ acting on a highest weight vector $v_{\lam} \in V^{\lam}.$ This is essentially contained in the following,

\lem{nonzero} Let $m \in \ZZ_{\geq 0}$ be a positive integer and $a=\alpha \otimes t^{-m} \in R_{-,re}$ be a root of depth $m$. Assume also that $\ome$ is regular and dominant. Then we have
\begin{enumerate}
  \item
  $\xi_{-m, \alpha} v_{\lam} \neq 0$
  \item
  If $h \in \Heis$ is of degree $-m$, then $h v_{\lam} \neq 0.$ Moreover, the family of vectors $\{ \xi_{-m, i} v_{\lam} \},$ $i=1, \ldots, l,$ are linearly independent.
  \item $\chi_{-a}(s) v_{\lam} = v_{\lam} + s \xi_{-m, \alpha} v_{\lam} +  \text{ terms of greater depth } $
 \end{enumerate}
\elem
\prf Parts (1) and (3) are clear. Let us prove (2). For any element of degree $n$ in $\Heis,$ say $h_n:= h \otimes t^n$ with $h \in \grt_o,$ we may find an element $h_{-n}:= h' \otimes t^{-n} \in \Heis$ with $h' \in \grt_o$ such that
$$
[ h_n, h_{-n} ] = \kappa \textsf{c}
$$
 with $\kappa \neq 0.$ Indeed, the this follows immediately from the commutation relation,
$$
[ h_n , h_{-n} ] = n (h | h') \textsf{c}
$$
and the non-degeneracy of $( \cdot | \cdot)$ when restricted to $\grt_o$ But we also have,
$$
[h_n, h_{-n}] v_{\lam} = h_{n} h_{-n} v_{\lam}
$$
 from which it must follow that $h_{-n} v_{\lam} \neq 0.$
\epr
\ssec{}{A reformulation}
In the theory of vertex operators, one encounters expressions of the following form for $i \in \{1, \ldots, l\}$ and
$p \in \ZZ_{\geq 0},$
$$
\mc{P}_i(s, p):= \exp( - \sum_{k > 0} \frac{ \xi_{pk,i} s^k }{k} ).
$$
It is easy to see that that $\mc{P}_i(s, p)$ is a well defined operator in $\Aut_\mathcal{K}(V^{\lambda})$ in the case when $\calK$ has
characteristic 0.
Moreover, we know from  \cite[Theorem 5.8, Remark (i)]{ga:la}
that if we define the element $\Lambda_k(\xi_i(p))$ via the relation,
$$
 \mc{P}_i(s, p) = \sum_{k\geq 0} \Lambda_k(\xi_i(p)) s^k
$$
then the elements $\Lambda_k(\xi_i(p))$ will actually lie in $U_{\ZZ}(\grg),$ our fixed $\ZZ$-form of the enveloping algebra. We then have the following,

\prop{} As elements of the corresponding completion of $U_{\ZZ}(\grg)$ we have an equality,
$$
\mc{P}_i(s, p) = h_i (1- s t^p).
$$
\eprop
\prf
We are going to check that the left hand side and the right hand side act in the same way in any irreducible integrable representation.
Let us first note that from \cite[Lemma 12.2]{gar:ihes}, it follows that $h_i(1 - s t^p)$ is in $U,$ and hence fixes any highest weight vector. From the definition of $\mc{P}_i(s, p),$ we see that it has the same property.Hence, by the analogue of Schur's Lemma \cite[Lemma 9.1]{gar:ihes}, in order to show that $\mc{P}_i(s, p)$ and $h_i(1- st^p)$ are equal, we need only to consider their action on the adjoint representation. But we may use the perfectness of $\grg$ as a Lie algebra and \cite[Lemma 8.11]{gar:ihes} to reduce to a computation in the adjoint representation of the loop algebra $\grg_o \otimes \mathcal{K}((t)).$ To complete the proof we know that, working in $\grg_o \otimes \mc{K}((t)),$ one has (in the notation above) that
$$
\Lambda_k(\xi_i(p)) x_{\beta} = (-1)^k \binom{ \langle \beta, \xi \rangle}{k} x_{\beta} \otimes t^{pk}
$$
where $x_{\beta} \in \mf{g_o}^{\beta}$ for $\beta \in R_o.$ By the usual binomial formula, one now has that
$$
\mc{P}_i(s, p)x_{\beta} = \sum_k (-1)^k \binom{ \langle \beta, \xi \rangle}{k} x_{\beta} \otimes t^{pk} s^k =
(1-t^ps)^{\langle \xi, \beta \rangle} x_{\beta}.
$$
 This agrees with the action of $h_i(1-st^p)$ on $x_{\beta}$ and so the lemma follows.
\epr

An entirely similar argument to the above shows the following,

\prop{}
The element $\mc{P}^-_i(s, p):= \exp( - \sum_{k > 0} \frac{ \xi_{-pk,i} s^k }{k} )$ defines an element of $\hatU^-$ and we have an equality of elements of $\hatU^-,$
$$
\mc{P}^-_i(s, p) = h_i(1- st^{-p}).
$$
\eprop

\cor{} Let $c_i \in k$ for $i=1, \ldots, l.$ Then there exists a positive integer $m_0$ so that for all $m \geq m_0$ we have we have an equality in $V[m],$  $$\prod_{i=1}^l h_i(1 - c_i t^{-m}) v_{\lam} = v_{\lam} + \sum_{i=1}^l c_i \xi_{m, i} v_{\lam} + \text{ terms of greater depth}.$$ Moreover, $\sum_{i=1}^l c_i \xi_{m, i} v_{\lam} \neq 0,$ provided not all $c_i=0$ for $i=1, \ldots, l.$
\ecor

The proof of the following theorem is parallel to that for $U$
described in \cite{gar:ihes}:
\th{imcoord}
Every element $u^- \in U^- \subset \hatU^-$ has a unique expression of the form \[ u^- = \prod_{m=0}^{\infty} u^-[m] = \cdots u^-[m] u^-[m-1] \cdots u^-[0] \]
(the product is considered in \emph{decreasing} order of m) where the $u^-[m]$ are of the form \refe{im} or \refe{im:0}. Furthermore, with respect to the map $\omega_m: \hatU^- \to U^-[m]$ we have
$$
\omega_m(u^-) = \omega_m(u^-[m] \; u^-[m-1] \; \cdots \; u^-[0]).
$$
\eth
Though there is no easy way to put coordinates on $U^-$, one can use the above theorem to define coordinates on $\hatU^-$.


\sec{part2}{Finiteness of the Gindikin-Karpelevich Sum, Part 2}
\noindent \emph{Notation:} For a constant $M,$ we write $M(x, y, z, ...)$ to indicate that the choice of $M$ depends only on $x, y, z, ....$
\vspace{0.15in}

\ssec{}{} Fix some highest weight module $V^{\lam}$ with primitive highest weight $v_{\lam}$.  Let us define a norm $|| \cdot ||$ on $V^{\lam}$ in teh following way.
For any $v\in V^{\lam}$ let us set
$$
\ord (v) = \min\limits_{n\in\ZZ}\pi^n v\in V^{\lam}_{\calO}.
$$
Then we set
$$
||v|| = q^{\ord (v)}.
$$
For an element $g \in G$ and a positive constant $C > 0$ we shall say that $g$ is
\emph{bounded by $C$} if $|| g v_{\lam} || < C.$
We say that $g$ is \emph{bounded by $C$ at depth $j$} if
$$
|| \sum_{\mu \in \calL_{\lam}(j) } v_{\mu} || < C.
$$
where
$$
gv_{\lam} = \sum_{\mu \in \calL_{\lam} } v_{\mu}
$$

We say that a family of elements $\Xi \subset G$ is bounded (bounded at depth $j$) if there exists some $C > 0$ such that every element of $\Xi$ is $C$-bounded (respectively, $C$-bounded at depth $j$).

 We shall need a finer notion of boundedness in the sequel.  Recall that each $u^-[j]$ can be
explicitly written  in terms  coordinates (see \refe{im} if $j >0$ and \refe{im:0} if $j=0.$)  We  say that an element $u^-[j]$ is \emph{componentwise bounded} by a positive number $C$ if in an expression \refe{im} or \refe{im:0} we have if $j > 0$
$$
 || s_{j, \alpha} || < C,  || \tilde{s}_{j, \alpha} || < C,  || c_{j, m} || < C \text{ where } \alpha \in R_o, \; m \geq 0
$$
or if $j=0$
$$
|| s_{0, \alpha} || < C \text{ where } \alpha \in R_o.
$$
A family of elements of the form $u^-[j]$ will be said to be componentwise bounded if there exists a constant $C$ such that all elements in this family are componentwise bounded by $C.$
\ssec{}{} The relation between between an element in $U^-$ being bounded and being componentwise bounded is partially explained by the following,

\lem{one:comp}
Suppose that $u^-[j]$ is bounded at depth $j$ by $C.$ Then there exists a constant $D=D(C)$ such that $u^-[j]$ is componentwise bounded by $D.$
\elem

\prf
We will consider two cases:

\ssec{}{Case 1: $j=0$}
Let us write
$$
u^-[0] = \prod\limits_{\alpha \in R_{o,-}} \chi_{\alpha}(s_{0, \alpha})
$$
where the product is ordered according to decreasing height from left to right.  Consider now the following statement for each $t=1, \ldots, \hite(\theta).$
\vspace{0.15in}

\noindent $\mc{H}(t):$ Suppose there exists $C >0$ such that for any
\eq{u:t}
u^-_t:= \prod_{\alpha \in R_{o, -} ,\; \hite(a) \geq t} \chi_{\alpha}(s_{0, \alpha}), \; \text{ } s_{0, \alpha} \in \mathcal{K} \end{equation}
that satisfies $$ || u_t v_{\lam} || < C $$ then there exists a constant $D= D(C) > 0$ (depending only on $C$) such that $|| s_{0, \alpha} || < D$ for each $\alpha \in R_{o, -}.$

\vspace{0.15in}

We shall show that $\mc{H}(t)$ is true for $t=1$ by decreasing induction on $t.$
For $t= \hite(\theta)$ this follows from \refl{nonzero}.
So we need to argue that that if $\mc{H}(t+1)$ is true, then also $\mc{H}(t)$ is true. So given $u^-_t$ as in
\refe{u:t} we may write
$$
u^-_t = u^-_{t+1} \chi_{\gamma_1}(s_{\gamma_1}) \cdots \chi_{\gamma_r}(s_{\gamma_r})
$$
where $s_{\gamma_j} \in \mathcal{K},$ $\hite(\gamma_j)=t$ for $j=1, \cdots r$ and $u^-_{t+1}$ is a product of elements corresponding to roots of height at least $t+1$.  We then have
$$
u^-_t v_{\lam} = v_{\lam} + s_i \xi_{\gamma_i} v_{\lam} + \text{ terms of lower order },
$$
By \refl{nonzero} each of the $\xi_{\gamma_i} v_{\lam} \neq 0,$ and furthermore they lie in different weight spaces.  Hence by the hypothesis of boundedness in $\mc{H}(t)$ we see that there exists a constant $D$ only depending on $C$ such that $|| s_{\gamma_j} || < D$ for $j=1, \ldots, r.$ Now, if we set
$$
\tilde{u}^-_{t+1}:= ( \chi_{\gamma_1}(s_{\gamma_1}) \cdots \chi_{\gamma_r}
 ) ^{-1} u^-_{t+1} \chi_{\gamma_1}(s_{\gamma_1}) \cdots \chi_{\gamma_r}, $$ then there exists a constant $\tilde{C}$ such that $|| \tilde{u}^-_{t+1} v_{\lam} || < \tilde{C}$ for all $u^-_t.$ On the other hand, the element $\tilde{u}^-_{t+1}$ is a product over roots of height at least $t+1$ and so we may apply inductively the hypothesis $\mc{H}(t+1)$ to the $\tilde{u}^-_{t+1}$ and conclude that it is componentwise bounded by a constant depending only on $C$. By as $\tilde{u}^-_{t+1}$ and $u^-_t$ are conjugates by a componentwise bounded expression, we see that the $u^-_t$ are also componentwise bounded.
\ssec{}{Case $j > 0$} Suppose we write
$$
u^-[j]=  \prod_{\alpha \in R_{o, +} } \chi_{\alpha}(t^{-j} s_{j, \alpha} ) \; \prod_{i=1}^l h_i(1+c_{i, j}t^{-j}) \; \prod_{\alpha \in R_{o, +}} \chi_{- \alpha} (t^{-j} \tilde{s}_{j, \alpha} ) \text{ with } s_{j, \alpha},\; \tilde{s}_{j, \alpha},\; c_{i, j} \in \mathcal{K}.
$$
Then by \refl{nonzero} we have
$$
u^-[j] v_{\lam} = v_{\lam} + \sum_{i=1}^l c_{i, j} \xi_{j, i} v_{\lam} + \sum_{\alpha \in R_{o, +} } s_{m, \alpha} \xi_{m, \alpha} v_{\lam} + \sum_{\alpha \in R_{o, +}} \tilde{s}_{m, -\alpha} \xi_{m, -\alpha} v_{\lam} + \text{ terms or lower depth} .
$$
Again by \refl{nonzero} we have that $ \xi_{j, i} v_{\lam}$ are linearly independent and so we may bound each of the coefficients $c_{i, j}.$ Furthermore, each of the vectors $\xi_{m, \alpha} v_{\lam}$ and $\xi_{m, -\alpha} v_{\lam}$ are non-zero and lie in different weight spaces, so we obtain a bound on the coefficients $s_{m, \alpha}$ and $\tilde{s}_{m, - \alpha}.$
\epr

An extension of \refl{one:comp} is given by the following,

\prop{uniform}
Let $m$ be a positive integer, and $C $ a positive constant. Suppose that $\mc{F} \subset U^-$ is a family such that \newline
\indent (a) $|| u^- v_{\lam} || < C$ for all $u^- \in \mc{F}$ \\
\indent (b) every $u^- \in \mc{F}$ may be written as a product $u^-[m] \cdots u^-[0].$  \newline  Then there exists a
$D=D(C) >0$ such that $u^-[j]$ is componentwise bounded by $D,$ for $j=0, \ldots, m.$
\eprop
\begin{proof} The proof will consist of a decreasing induction on $j$ for the following statement denoted by $\mc{P}(j)$:

\vspace{0.1in}

\noindent $\mc{P}(j)$: Suppose that $u^-[m] u^-[m-1] \cdots u^-[j]$ is bounded by $C.$ Then there exists $D=D(C)$ such that each $u^-[k]$ for $k=j, \cdots, m$ is componentwise bounded by $D.$

\vspace{0.1in}

The statement $\mc{P}(m)$ follows from  \refl{one:comp}.  So let us assume that $\mc{P}(j+1)$ is true, and let us then argue that $\mc{P}(j)$ then follows. Let us write $$ u^-:= u^-[m] \cdots u^-[j+1] u^-[j] v_{\lam} = v_{\lam} + \text{ terms of depth $j$ } + \text{ terms of higher depth }. $$ We then have that $$ (u^- - u^-[j])v_{\lam} =   \text{ terms of depth  $\geq  j+1$}. $$ As $u^-$ belongs to a bounded family, we must have that the $u^-[j]$ is bounded at depth $j.$ Hence it is bounded and componentwise bounded by \refl{one:comp}. Now consider $$ \tilde{u}^-:= u^-[j]^{-1} u^-[m] \cdots u^-[j] = \tilde{u}^-[m] \cdots \tilde{u}^-[j+1] $$ for some elements $\tilde{u}[m], \ldots, \tilde{u}[j+1].$  The expression $\tilde{u}$ is bounded by some constant $D=D(C)$ and so applying statement $\mc{P}(j+1)$ we conclude that $\tilde{u}$ is componentwise bounded by $E=E(D)=E(C).$ But then there exists a constant $F=F(E)=F(D)=F(C)$ such that $$ u^- = u^-[j] \tilde{u}^- u[j]^{-1} $$ is componentwise bounded (since $u^-[j]$ was componentwise bounded).

\end{proof}
In the future we are going to need the following result:
\prop{AB}
Let $n$ be a positive integer and let $C>0$.
There exists $r=r(C, n) \geq 0$ such that for any  $A, B\in GL(n,\calK)$ such that
\begin{enumerate}
  \item The entries of $A,B$ are bounded by $C$
  \item $A-B\equiv 0\; (\md \pi^r)$,
\end{enumerate}
Then $AB^{-1}\in GL(n,\calO)$.
\eprop
The proof is easy and it is given in the Appendix.
We now proceed to the proof of \reft{main-fin}.

\ssec{}{Step 1}Let $\lam^{\vee}$ be a regular dominant weight. Then the natural map $\zeta: U^-\to V^{\lam}$ sending
every $u^-$ to $u^-(v_{\lam})$ is injective. We claim that the same is true for
the map $\zeta_m:U^-(w,\mu^{\vee})\to V^{\lam}[m]$ (which is equal to the composition of $\zeta$ with the projection
 $V^{\lam}\to V^{\lam}[m]$) for sufficiently large $m$.
Indeed, since $U(w,\mu^{\vee})\subset BwB$, it follows that
the image of $\zeta$ lies in $V_w^\lam$. Since the projection $V_w^\lam\to V^{\lam}[m]$ is injective for large $m$, the statement
follows.

Moreover, we claim that if $m$ is sufficiently large, then both $\zeta_m$ and the map
$U^-\supset U^-(w,\mu^{\vee})\cdot U^-(w,\mu^{\vee})\to V^{\lam}[m]$ sending $(u_1^-,u_2^-)$ to
$u_2^- (u_1^-)^{-1} v_{\lam} \;(\md V^{\lam}(m) ) $ are injective. Indeed, if $u_1^-, u_2^- \in U^-(w, \mu)$ then the product
$$
u_2^- (u_1^-)^{-1} \in BwBw^{-1} B \subset \bigcup_{w' \in \Omega} B w' B
$$
where $\Omega$ is a finite set. Then any $m\geq \max_{w' \in \Omega} m_{w'}$ will satisfy the second injectivity requirement.

The second step is the following simple lemma:
\lem{Step 2}The set $U(w,\mu^{\vee})$ is bounded. Equivalently, there exists  $i>0$ such that
$\zeta(U^-(w,\mu^{\vee}))\subset \pi^{-i}V^{\lam}_{\calO}$.
\elem
\prf Take $i=\la \mu^{\vee},\lam\ra$. Then  for every $g\in K\pi^{-\mu^{\vee}}U$ we have
$$
\zeta(g)\subset \pi^{-i}V^{\lam}_{\calO}.
$$
Hence \refl{Step 2} follows since
$U^-(w,\mu^{\vee})\subset K\pi^{-\mu^{\vee}}U$.
\epr

\ssec{}{Step 3}Let $m$ be greater than the depth of $w(\lam)$.
Then for every $u^-\in U(w,\mu^{\vee})$  in the decomposition
\eq{bd}
 u^- v_{\lam} = \sum_{\mu \in P_{\lam} } v_{\mu} , v_{\mu} \in V^{\lam}(\mu),
\end{equation}
we must have $v_{\mu} = 0$ for all $\mu$ of $\dep(\mu) > m.$

Using \reft{imcoord} let us write $u^- = \prod_{j=0}^{\infty} u^-[j].$
Then for every $m$ as above
we see immediately that $\omega_m(u^-[j])=1$ if $j \geq m.$  Hence for sufficiently large $m$ (independent of $u^-$ we have
\eq{ufin}
u^- = u^-[m] u^-[m-1] \cdots u^-[0]
\end{equation}
\ssec{}{Step 4}Let us choose $m$ to satisfy the conditions of Step 1 and Step 3.
From now on let us set $V=V^{\lam}_w, V_{\calO}=V^{\lam}_{w,\calO}, v=v_{\lam}$. We claim
that there exists $c>0$ such that for any $u^-\in U^-(w,\mu^{\vee})$ we have
$$
u^-(V_{\calO})\subset \pi^{-c} V_{\calO}.
$$
In other words, we claim that if we choose an $\calO$-basis for $V_{\calO}$ then the image of $U(w,\mu^{\vee})$ in $GL(V)$ consists of elements whose entries are bounded by some constant $C=q^c$. This immediately follows from \refp{uniform}.

\ssec{}{Step 5}
We want to show that
there exists a finite set $\mc{F} \subset U^-(w, \mu^{\vee})$ such that for any $u^- \in U^-(w, \mu^{\vee})$ there exists $u^-_{\O} \in U^-_{\O}$ and $u^-_f \in \mc{F}$ such that
$$
u^-=u^-_{\O} u^-_f.
$$
For every positive integer $l \geq 0$ we construct a set $\mc{F}_l$ as follows. We have already constructed an embedding $\xi$
of $U^-(w,\mu^{\vee})$ into the group of  unipotent lower triangular matrices whose entries are uniformly bounded by some constant $C.$
For a given $l \geq 0$ there are only finitely many such matrices $\md \pi^l.$  Let $\mc{F}_l\subset U^-(w,\mu^{\vee})$ be a set of representatives of $\xi(U^-(w, \mu))$ mod $\pi^l$. In other words, for every element $u^- \in \mc{F}_l,$ we may write
\eq{princ:part}
\omega_l(u^-) = A_0 + \epsilon \; \text{ where } \epsilon \equiv 0 \; (\md \pi^l), A_0 \in \mc{A}_l.
\end{equation}
We claim that if $l$ is sufficiently large, then the set $\calF_l$ satisfies our requirements.
First we choose $l$ satisfying the condition of \refp{AB}.  Then for every $u^- \in U^-(w, \mu),$ there exists $u^-_f \in \mc{F}_l$ the such that $\xi((u^-))^{-1}\xi( u^-_f)$ has integral entries.
So, we know that $ (u^-)^{-1} u^-_f v_{\lam} \in V^{\lam}_{\O}[m]$ and hence also $(u^-)^{-1} u^-_f v_{\lam} \in V^{\lam}_{\O}$ by
the second injectivity requirement from Step 1. Hence, it is enough to prove the following
\lem{int-int}
Let $\lam$ be a regular dominant weight. Assume that $u^- v_{\lam}\in V^{\lam}_{\calO}$ for some $u^-\in U^-$. Then
$u^-\in U^-_{\calO}$.
\elem
\prf
For simplicity let us assume that $\lam=\rho$ (the general case is similar). 
We know that $u^-\in K\pi^{-\mu^{\vee}}U$ for some $\mu^{\vee}\in Q^{\vee}_+$. 
The fact that $u^-v_{\lam}\in V^{\lam}_{\calO}$ implies that 
$\mu^{\vee}=0$. Hence $u^-\in U^-_{\calO}$ by \refl{drin}.
\epr
\sec{final}{Proof of \reft{finiteness} and \reft{gk-aff}}
\ssec{}{Proof of \reft{finiteness}(1) and \reft{finiteness}(3)}
Let $\lam\in \Lam^+$. Then any $g\in K\pi^{\lam^{\vee}}K$ satisfies
\eq{kpk}
g^{-1}(V^{\lam}_{\calO})\subset \pi^{-\la\lam^{\vee},\lam\ra} V^{\lam}_{\calO}.
\end{equation}
Indeed, the condition \refe{kpk} is clearly $K$-bi-invariant and it is trivially satisfied by
$g=\pi^{\lam^{\vee}}$.

On the other hand, any $g\in K\pi^{\mu^{\vee}} U^-$ satisfies
\eq{kpum}
\eta_{\lam}(g^{-1}(V^{\lam}_{\calO}))\subset \pi^{-\la\mu^{\vee},\lam\ra}\calO,
\end{equation}
where $\eta_{\lam}:V^{\lam}\to\calK$ is the projection to the highest weight line (normalized by the condition
that $\eta_{\lam}(v_{\lam})=1$). This is true because the set of all $g$ that satisfy \refe{kpum} is clearly
invariant under $U^-$-action on the right and the $K$-action on the left and it is satisfied by $g=\pi^{\mu^{\vee}}$.

Thus we see that if
$$
K\pi^{\lam^{\vee}}K\cap K\pi^{\mu^{\vee}} U^-\neq\emptyset ,
$$
then \refe{kpk} and \refe{kpum} imply that $\la\lam^{\vee},\lam\ra\geq \la\mu^{\vee},\lam\ra$ for every $\lam\in \Lam^+$. Hence
$\lam^{\vee}\geq \mu^{\vee}$ which proves \reft{finiteness}(3).

Similarly, if $g\in K\pi^{\mu^{\vee}}U$ then
\eq{kpu}
g^{-1}(v_{\lam})\in \pi^{-\la\mu^{\vee}, \lam\ra} V^{\lam}_{\calO}.
\end{equation}
Hence \refe{kpum} and \refe{kpu} imply that if
$$
K\pi^{\lam^{\vee}}U^-\cap K\pi^{\mu^{\vee}} U\neq\emptyset ,
$$
then $\la\lam^{\vee},\lam\ra\geq \la\mu^{\vee},\lam\ra$ for every $\lam\in \Lam^+$. Hence
$\lam^{\vee}\geq \mu^{\vee}$. This proves the second assertion of \reft{finiteness}(1) and we already know the first
assertion.
\ssec{}{Proof of \reft{finiteness}(2)}
The second assertion of \reft{finiteness}(2) follows immediately from \refe{kpk} and \refe{kpu}.
Now we need to prove that every $K\pi^{\lam^{\vee}}K\cap K\pi^{\mu^{\vee}}U$ is finite.
But we can write
\eq{triple}
K\pi^{\lam^{\vee}}K\cap K\pi^{\mu^{\vee}}U=\bigcup\limits_{\lam^{\vee}\geq\nu^{\vee}\geq \mu^{\vee}}
K\pi^{\lam^{\vee}}K\cap K\pi^{\nu^{\vee}}U^-\cap K\pi^{\mu^{\vee}}U.
\end{equation}
Each $K\pi^{\lam^{\vee}}K\cap K\pi^{\nu^{\vee}}U^-\cap K\pi^{\mu^{\vee}}U$ is finite since it is a subset in the finite
set $K\pi^{\nu^{\vee}}U^-\cap K\pi^{\mu^{\vee}}U$. Since there are finitely many $\nu^{\vee}$ such that
$\lam^{\vee}\geq\nu^{\vee}\geq \mu^{\vee}$, it follows that the right hand side of \refe{triple} is finite.
\ssec{}{Proof of \reft{finiteness}(4)}This statement is well-known in the finite case and the proof in the affine case
is essentially similar. We include it here for completeness.

We need to prove the following
\prop{reform}
Let us fix $\mu^{\vee}\in Q^{\vee}_+$. Then
 for sufficiently dominant $\lam^{\vee}$ we have

(1)  $K\pi^{\lambda^{\vee}} K \cap K \pi^{\lambda^{\vee} - \mu^{\vee}} U \subset K \pi^{\lambda^{\vee}} U^- $

(2)$ K \pi^{\lambda^{\vee}} U^-\cap K \pi^{\lambda^{\vee} - \mu^{\vee}} U \subset K\pi^{\lambda^{\vee}} K $
\eprop

\prf
If $\lambda^{\vee}$ is dominant, then $K \pi^{\lambda^{\vee}} I \subset K \pi^{\lambda^{\vee}} U^-$.
So in order to prove (1) it is enough to prove the following
\lem{claim:bg}
For sufficiently dominant $\lambda^{\vee}$  we have
 $$
 K \pi^{\lambda^{\vee}} K \cap K \pi^{\lambda^{\vee} - \mu^{\vee}} U \subset K \pi^{\lambda^{\vee}} I
 $$
\elem
\prf
Suppose we have an element $ x \in K \pi^{\lambda^{\vee}} K \cap K \pi^{\lambda^{\vee} - \mu^{\vee}} U.$ So we may write
$$
x = k_1 \pi^{\lambda^{\vee}} k_2 = k_3 \pi^{\lambda^{\vee} - \mu^{\vee}} u
$$
where $k_i \in K, u \in U.$ This then implies that there exists $k_4 \in K$ such that
\eq{claim:bg:1}  \pi^{\lambda^{\vee}} k_2 = k_4 \pi^{\lambda^{\vee} - \mu^{\vee}} u.
\end{equation}
We need to show if $\lambda^{\vee}$ is sufficiently dominant,
then $k_2 \in I.$  In other words, we need to show that if $k_2 \in I w I$ for $ w \in W$ then $w=1.$

Let us choose a dominant weight $\lam$ and let $V^{\lam}, v_{\lam}$ be as before.
Then for  sufficiently dominant $\lambda^{\vee}$  the following condition is satisfied:
\eq{suff}
\text{ if } w \neq 1, \text{ and } w \lambda^{\vee} = \lambda^{\vee} - \beta^{\vee}, \beta^{\vee} \in Q_+, \text{ then } \la \beta^{\vee}, \lambda \ra > \la \mu^{\vee}, \lambda \ra.
\end{equation}
Rewriting \refe{claim:bg:1} under the assumption that $k_2 \in I w I,$ we have
$$
\pi^{\lambda^{\vee}} i_1 w i_2 \in K \pi^{\lambda^{\vee} - \mu^{\vee}} U ,\; \text{ for } i_1, i_2 \in I
$$
Since $I= U^+_{\O}U^-_{\pi}T_{\O}$ and $\pi^{\lambda^{\vee}}U^+_{\O}T_{\O} \pi^{-\lambda^{\vee}} \subset U^+_{\O}T_{\O}$ and $w U^-_{\pi} w^{-1} \subset I$ we have that
$$
\pi^{\lambda^{\vee}} I w I = \pi^{\lambda^{\vee}} U^+_{\O} U^-_{\pi} T_{\O} w I \subset K \pi^{\lambda^{\vee}} w U^-_{\pi} U^+_{\O}
$$
So if \refe{claim:bg:1} holds, then we can conclude that there exists $u^- \in U^-_{\pi},$ $k \in K$ and $u \in U$ such that \eq{bg:claim:2}
\pi^{\lambda^{\vee}} w u^- = k \pi^{\lambda^{\vee} - \mu^{\vee}} u
\end{equation}
Now, apply \refe{bg:claim:2} to the highest weight vector $v_{\lambda}:$ \eq{bg:claim:3}
\pi^{\lambda^{\vee}} w u^- v_{\lambda} = k \pi^{\lambda^{\vee} - \mu^{\vee}} u v_{\lambda}
\end{equation}
Consider the left hand side of
\refe{bg:claim:3} we obtain,
 $$
 ||  \pi^{\lambda^{\vee}} w u^- v_{\lambda} || \geq || \pi^{\lambda^{\vee}} w v_{\lambda} || = q^{- \langle \lambda^{\vee}, w \lambda \rangle }.
 $$
From the right hand side of \refe{bg:claim:2} we obtain,
$$
 ||k \pi^{\lambda^{\vee} - \mu^{\vee}} u v_{\lambda} || = q^{- \langle \lambda^{\vee} - \mu^{\vee}, \lambda \rangle }.
$$
Writing $w^{-1} \lambda^{\vee} = \lambda^{\vee} - \beta^{\vee}$ for $\beta^{\vee} \in Q_+,$ we have that
$$
q^{- \langle \lambda^{\vee}, w \lambda \rangle} = q^{- \langle w^{-1} \lambda^{\vee}, \lambda \rangle} = q^{- \langle \lambda^{\vee} - \beta^{\vee}, \lambda \rangle} \leq q^{- \langle \lambda^{\vee} - \mu^{\vee}, \lambda \rangle }.
$$
This implies that
$$
 q^{\langle \beta^{\vee}, \lambda \rangle} \leq q^{\langle \mu^{\vee}, \lambda \rangle}
$$
which contradicts the fact that we have chosen $\lambda^{\vee}$ to satisfy \refe{suff}.

\epr
\ssec{}{Prof of \refp{reform}(2)}
There exists a finite set $\Omega \subset U^-$ such that
$$
K \pi^{\lambda^{\vee} - \mu^{\vee}} U \cap K \pi^{\lambda^{\vee}} U^- \subset \bigcup_{u^- \in \Omega} K \pi^{\lambda^{\vee}} u^-. $$
If all such $u^- \in \Omega$ actually lie in $K$ we are done. Let $\Omega_0 \subset \Omega$ be the (finite) subset of $u^- \in \Omega$ such that $u^- \notin K.$ Let $u^- \in \Omega_0.$ Since $u^- \notin K$  there exists $w\in W$ with $l(w) > 1$ such that $u^- \in U^-_{w},$ where $U^-_w$ is as in \refs{part1}.
Let us write
$$
u^- v_{\rho} = v_{\rho} + \sum_{\gamma \in Q_+} v_{\rho - \gamma}.
$$
Then there are only finitely many $\gamma$ which appear in such expressions as $u^- \in \Omega.$
By  \refl{drin}, we have that $|| u^- v_{\rho} || > 0,$ so there exists $\gamma \neq 0$ such that $v_{\rho - \gamma} \notin V^{\rho}_{\O}.$

By the hypothesis, we have an expression of the form
$$
\pi^{\lambda^{\vee}} u^- = k  \pi^{\lambda^{\vee} - \mu^{\vee}} u^+ .
$$
Applying the right hand side to $v_{\rho}$ we find that
$$
|| k \pi^{\lambda^{\vee} - \mu^{\vee}} u^+ v_{\rho} || = q^{- \la \lambda^{\vee} - \mu^{\vee}, \rho \ra }
$$
whereas when we apply the same expression to the left hand side, we obtain that
$$
|| \pi^{\lambda^{\vee}} u^- v_{\rho} || \geq || \pi^{\lambda^{\vee}} v_{\rho - \gamma} || \geq q^{ - \la \lambda^{\vee}, \rho - \gamma \rangle} .
$$
By choosing $\lambda^{\vee}$ sufficiently dominant so that $\la \lambda^{\vee} , \gamma \ra > \la \rho, \mu^{\vee} \ra$ for the finitely many $\gamma$ which can occur, we obtain a contradiction.

\epr

\ssec{}{Proof of \reft{gk-aff}}
Let us prove \reft{gk-aff} assuming \reft{aff-mac}. It follows from \reft{finiteness}(4)
that the coefficient of some $e^{\gam^{\vee}}$ in $\bfc_{\grg}$ is equal to the coefficient
of $e^{\lam^{\vee}-\gam^{\vee}}$ in $\calS(\lam^{\vee})$ divided by $q^{\la \rho,\lam^{\vee}\ra}$
for sufficiently large $\lam^{\vee}$.
However, it is clearly equal to the right hand side of \refe{gk-aff} (since only the term with $w=1$
in the formula \refe{hlam}
can contribute to the coefficient $e^{\lam^{\vee}-\gam^{\vee}}$ for sufficiently large $\lam^{\vee}$; also
if $\lam^{\vee}$ is sufficiently large, then it is also regular and thus $W_{\lam^{\vee}}$ is trivial).

\sec{app}{Appendix: Proof of \refp{AB}}

\ssec{}{}  Let $A=(a_{ij}) $ be an $n \times n$ lower triangular, unipotent matrix with coefficients bounded by $|| a_{i, j} || < C.$  Then $A^{-1}$ is also a bounded matrix whose entries are bounded by some $C'=C'(C, n).$  For $k \geq 0,$ assume we have an expression
\eq{A0}  A = A_0 + \epsilon
\end{equation}
where $A_0$ is lower triangular unipotent and $\epsilon$ are $n \times n$ matrices such that the coefficients of $\epsilon$ all lie in $\pi^{k} \O.$

\lem{ep-del}
Let $A$ and $A_0$ be $n \times n$ unipotent lower triangular matrices whose coefficients have norm bounded by some number $C.$  Then given any $m \in \ZZ_{\geq 0}$, there exists $l_0=l_0(C, n, m)$ such that if $A=A_0 + \epsilon$ with $\epsilon \equiv 0 (\md \pi^{l_0})$ then
$$
A^{-1} = A_0^{-1} + \delta
$$
with $\delta \equiv 0 (\md \pi^m).$
\elem
\begin{proof} We first note that the following two facts:

(1) Given $p \geq 0$, there there exists a positive integer $d_0=d_0(C, n, p)$ such that if $\epsilon \equiv 0 \; (\md \pi^{d_0} ),$ then $\epsilon A^{-1} \equiv 0 \;(\md \pi^p).$

\renewcommand{\id}{\mathbb{I}_n}

(2) The Lemma is true when $A_0=\id,$ the $n \times n$ identity matrix.

Indeed, (1) follows from the fact that if the coefficients of $A$ are bounded by $C$ then those of $A^{-1}$ must be bounded by some $C'.$ Also, (2) follows from the identity of matrices, $$ (\id + \epsilon)^{-1} = \id - \epsilon + \epsilon^2 + \cdots. $$

Now if we write
\eq{A:A0}
A= A_0 + \epsilon,
\end{equation}
then we have
$$
A^{-1} = (\id + A_0^{-1} \epsilon)^{-1} A_0^{-1}.
$$
Choosing $\epsilon$ sufficiently small (in the $\pi$-adic topology), then thanks to (1) and (2) we can assume that
$$
(\id + \epsilon A_0^{-1})^{-1} = \id + \tilde{\delta}
$$
where $\tilde{\delta}$ is sufficiently small.  Hence,
$$
A^{-1} = A_0^{-1}(\id + \tilde{\delta}) = A_0^{-1} + A_0^{-1} \tilde{\delta}
$$
and by stipulating that $\tilde{\delta}$ is sufficiently small we may assume that
$$
\delta = A_0^{-1} \delta
$$
is arbitrarily small.  In sum, choosing $\epsilon$ sufficiently small in \refe{A:A0} we have $A^{-1} = A_0^{-1} + \delta$ with $\delta$ arbitrarily small (uniformly for all $A, A_0$ lower triangular unipotent matrices with coefficients bounded by $C.$)
\end{proof}


\ssec{}{Proof of \refp{AB}}
Let us write
$$
A = X_0 + \epsilon_1,  B= X_0 + \epsilon_2,  \epsilon_i \equiv 0 \;(\md \pi^r) \; \text{ for } i=1, 2.
$$

Given a $p \geq 0,$ we may use \refl{ep-del} to choose $r=r(p, C, n)$ sufficiently large so that
$$
B^{-1} = X_0^{-1} + \delta \text{ for } \delta \equiv 0 (\md \pi^p).
$$
For this choice of $r$ we have
$$
AB^{-1} = (X_0 + \epsilon_1) (X_0^{-1} + \delta) = 1 + X_0 \delta + \epsilon_1 X_0^{-1} + \epsilon_1 \delta.
$$
Pick $p$ sufficiently large so that $X_0 \delta$ (and also trivially $\epsilon_1 \delta$) is integral.  Increasing the value of $r$ if necessary, we can also arrange that $\epsilon_1 X_0^{-1}$ is also always integral.  The proposition follows.

 \end{document}